\documentstyle{amsppt}
\magnification=1200
\hsize=150truemm
\vsize=224.4truemm
\hoffset=4.8truemm
\voffset=12truemm

\NoRunningHeads

\define\C{{\bold C}}
 
\let\thm\proclaim
\let\fthm\endproclaim
 

\define\Ga{\Gamma}
 \define\Gn{\Gamma_n}
\define\Si{\Sigma}
\define\Sn{\Sigma_n}
\define\Bn{\partial \Sigma_n}
 
\newcount\tagno
\newcount\secno
\newcount\subsecno
\newcount\stno
\global\subsecno=1
\global\tagno=0
\define\ntag{\global\advance\tagno by 1\tag{\the\tagno}}

\define\sta{\ 
{\the\secno}.\the\stno
\global\advance\stno by 1}

\define\stas{\the\stno
\global\advance\stno by 1}

\define\sect{\global\advance\secno by 1
\global\subsecno=1\global\stno=1\
{\the\secno}. }

\def\nom#1{\edef#1{{\the\secno}.\the\stno}}
\def\inom#1{\edef#1{\the\stno}}
\def\eqnom#1{\edef#1{(\the\tagno)}}

\newcount\refno
\global\refno=0

\def\nextref#1{
      \global\advance\refno by 1
      \xdef#1{\the\refno}}

\def\bref {\ref\global\advance\refno by 1\key{\the\refno}}

\nextref\THE
\nextref\NEV
\nextref\REY

\topmatter

\title 
Sur la th\'eorie d'Ahlfors des surfaces
 \endtitle

\author  Julien Duval\footnote""{Laboratoire de Math\'ematiques, Universit\'e Paris-Sud, 91405 Orsay cedex, France \newline julien.duval\@math.u-psud.fr\newline}
\footnote""{Mots-cl\'es :  covering surfaces, Riemann-Hurewitz formula, value distribution theory \newline Class.AMS : 30C25,30D35.\newline}
\endauthor

\abstract\nofrills{\smc R\'esum\'e. }\ \ On revisite la th\'eorie d'Ahlfors de recouvrement des surfaces \`a l'aide du th\'eor\`eme de Stokes.
 
\endabstract

\endtopmatter 

\document
 
\thm{Contexte}\fthm On pr\'esente la th\'eorie d'Ahlfors (voir par exemple [\THE], [\NEV]) de mani\`ere asymptotique. C'est sa forme utile et les formules sont plus simples. Pour des rappels sur les surfaces de Riemann on renvoie \`a [\REY].

Soit $f_n:\Sn \to \Si$ une suite d'applications holomorphes non constantes entre deux surfaces de Riemann compactes connexes $\Sn$ (avec bord) et $\Si$ (sans bord). On munit $\Si$ d'une m\'etrique d'aire totale 1 qui induit via $f_n$ une pseudom\'etrique sur $\Sn$. Dans toute la suite on suppose que la longueur $l_n$ du bord de $\Sn$ devient n\'egligeable devant son aire $a_n$. 

\thm{Hypoth\`ese} $l_n=o(a_n)$. \fthm

\noindent
On note $\lesssim$ ou $\sim$ une in\'egalit\'e ou \'egalit\'e \`a $o(a_n)$ pr\`es.

En cons\'equence, si $\theta$ est une 2-forme sur $\Si$ d'int\'egrale 1 alors $\int_{\Sn}f_n^*\theta\sim a_n$. En effet $\theta$ est cohomologue \`a la forme d'aire $\omega$, donc $\theta-\omega=d\alpha$ et par le th\'eor\`eme de Stokes $\int_{\Sn}f_n^*\theta-a_n=\int_{\Bn}f_n^*\alpha$ qui est contr\^ol\'ee par $l_n$. Une autre mani\`ere d'\'ecrire ceci est $\int_{\Si}d_n\theta\sim a_n $ o\`u $d_n(p)$ est le nombre de points dans $f_n^{-1}(p)$.

Dans ce contexte la th\'eorie d'Ahlfors dit que $f_n$ se comporte presque comme un rev\^etement ramifi\'e de degr\'e $a_n$. On a une in\'egalit\'e
de Riemann-Hurwitz sur les caract\'eristiques d'Euler.
 
\thm{In\'egalit\'e asymptotique}      $ \chi(\Sn) \lesssim a_n\chi (\Si)$. \fthm

Une version relative de cette in\'egalit\'e donne le th\'eor\` eme des \^\i les, qui est le pendant g\'eom\'etrique de la th\'eorie de Nevanlinna de distribution des valeurs. On suppose ici que $\Sn$ est un disque et $\Si$ la sph\`ere de Riemann. Fixons trois disques topologiques disjoints dans $\Si$ et appelons {\it \^\i le} toute composante connexe de la pr\'eimage par $f_n$ d'un de ces disques sur laquelle $f_n$ est propre. Soit $i_n$ le nombre total d'\^\i les au-dessus des disques.

\thm{Th\'eor\`eme des \^\i les} $ i_n\gtrsim a_n. $ \fthm

En particulier une telle suite ne peut \'eviter trois points distincts. Par un argument longueur-aire (voir plus bas) Ahlfors en d\'eduit le th\'eor\`eme de Picard, l'absence d'application enti\`ere non constante \`a valeurs dans la sph\`ere de Riemann priv\'ee de trois points. De la m\^eme mani\`ere l'in\'egalit\'e asymptotique donne l'absence d'application enti\`ere non constante \`a valeurs dans $\Si$ si $\chi(\Si)<0$.
 
\null
On va voir que ces deux r\'esultats proviennent du fait que $f_n$ se comporte vraiment comme un rev\^etement non ramifi\'e de degr\'e $a_n$ au-dessus d'un graphe $\Ga$ dans $\Si$, quitte \`a le perturber un peu.
\thm{\'Egalit\'e asymptotique} $ \chi(f_n^{-1}(\Ga))\sim a_n\chi(\Ga). $ \fthm

Ici $\Ga$ est une r\'eunion finie d'arcs $\gamma$ se coupant transversalement. La perturbation (d\'ependant de n) s'obtient en d\'epla\c cant un peu les arcs parall\`element \`a eux-m\^emes. La pr\'eimage $f_n^{-1}(\gamma)$ d'un arc \'evitant les valeurs critiques de $f_n$ se scinde en {\it bons arcs} (ceux se projetant hom\'eomorphiquement via $f_n$ sur $\gamma$) et en {\it mauvais arcs} (ceux rencontrant le bord $\Bn$). L'\'egalit\'e asymptotique r\'esulte directement du fait suivant.

\thm{Fait} Apr\`es perturbation $f_n^{-1}(\gamma)$ contient $\sim a_n$ bons arcs et $o(a_n)$ mauvais arcs.\fthm
 
Pour ceci param\'etrons un voisinage de $\gamma$ par un rectangle fin horizontal $R$ du plan. Les perturbations de $\gamma$ apparaissent comme des segments horizontaux $\gamma_t=\pi^{-1}(t)$ o\`u $\pi $ est la projection sur l'axe vertical. Le nombre de mauvais arcs au-dessus de $\gamma_t$ est estim\'e par le nombre de points de $f_n(\Bn)$ dans $\gamma_t$. Or $ \int\text{card}(f_n(\Bn) \cap \gamma_t)dt=\text{long}(\pi(f_n(\Bn)\cap R))$ o\`u la longueur est compt\'ee avec multiplicit\'e. Cette int\'egrale est donc contr\^ol\'ee par $l_n$ et on peut trouver une perturbation $\gamma_{t_n}$ avec $o(a_n)$ mauvais arcs. Par le lemme de Sard on peut m\^eme supposer $\gamma_{t_n}$ transverse \`a $f_n(\Bn)$.

\null
Regardons les bons arcs au-dessus de $\gamma_{t_n}$. Rappelons que $d_n(p)=\text{card}(f_n^{-1}(p))$. Soit $\beta$ une 1-forme  \`a support compact dans $R$ ne d\'ependant que de la variable horizontale et d'int\'egrale 1 sur $\gamma_t$. On va voir que $\int_{\gamma_{t_n}}d_n\beta\sim a_n$. Cela suffit car ainsi $d_n\gtrsim a_n$ en un point de $\gamma_{t_n}$. Les mauvais arcs \'etant peu nombreux, on a donc $\gtrsim a_n$ bons arcs au-dessus de $\gamma_{t_n}$, puis $\sim a_n$ par l'estim\'ee int\'egrale. Montrons-la. Soit $\delta(t)$ une fonction d'int\'egrale 1 \`a support compact et $\delta_{\epsilon}=\frac 1 {\epsilon} \delta\circ h_{\epsilon}$ l'approximation de l'unit\'e correspondante en $t_n$. On a $\beta \wedge\delta_{\epsilon}dt-\omega=d\alpha_{\epsilon}$  o\`u $\alpha_{\epsilon}$ est une 1-forme {\it uniform\'ement born\'ee}. Par exemple $\alpha_{\epsilon}=\alpha +(\phi-\phi\circ h_{\epsilon})\beta$ o\`u $\phi$ est une primitive de $\delta$ et $d\alpha=\beta \wedge \delta dt-\omega$. Donc $\int (\int_{\gamma_t}d_n\beta)\delta_{\epsilon}dt-a_n$ est contr\^ol\'ee uniform\'ement par $l_n$. Comme $t\to \int_{\gamma_t}d_n\beta$ est continue en $t_n$ (par transversalit\'e de $\gamma_{t_n}$ au bord), l'estim\'ee s'obtient par passage \`a la limite en $\epsilon$.  

\thm{Pour le th\'eor\`eme des \^ \i les}\fthm Voici comment le d\'eduire. On place un graphe $\Ga$ en forme de figure huit dans la sph\`ere de Riemann de sorte que chacune des composantes connexes du compl\'ementaire de $\Ga$ contienne un des disques topologiques. Par l'\'egalit\'e asymptotique $\chi(f_n^{-1}(\Ga))\sim -a_n$ car $\chi(\Ga)=-1$. Soit $\Gn$ le graphe obtenu en supprimant les mauvais arcs de $f_n^{-1}(\Ga)$. On a encore $\chi(\Gn)\sim-a_n$. Le compl\'ementaire de $\Gn$ dans le disque $\Sn$ consiste en une composante connexe $C_0$ touchant le bord et d'un certain nombre de composantes connexes $C$ \`a l'int\'erieur de $\Sn$. N\'ecessairement $\chi(C_0)\leq0$ et $\chi(C)\leq1$ car $C_0$ a au moins deux composantes de bord et $C$ une. Par ailleurs, par construction chaque composante $C$ contient une \^\i le. En effet $f_n(\partial C) \subset \Ga$, donc $f_n(C)$ couvre une composante du compl\'ementaire de $\Ga$ par le principe du maximum. L'\'egalit\'e en caract\'eristique d'Euler
$\chi(\Sn)=\chi(C_0)+\chi(\Gn)+\sum\chi(C)$
donne $1\lesssim -a_n +i_n$.

\thm{Pour l'in\'egalit\'e asymptotique}\fthm
Elle suit le m\^eme sch\'ema. On prend pour $\Ga$ une dissection de $\Si$, i.e. un bouquet de $2g$ cercles (si $g$ est le genre de $\Si$) dont le compl\'ementaire est un disque topologique. On a $\chi(\Si)=\chi(\Ga)+1$. On construit $\Gn$, les composantes $C_0$ touchant $\Bn$ et $C$ comme plus haut. Par l'\'egalit\'e asymptotique $\chi(\Gn)\sim a_n\chi(\Ga)$.
Par ailleurs le nombre de composantes $C$ est au plus $a_n$. En effet une telle composante se projette proprement sur le disque donc $\int_Cf_n^*\omega =
\text {deg}(f_n\vert_C)\geq1$. L'\'egalit\'e en caract\'eristique d'Euler donne $\chi(\Sn)\leq \chi(\Gn)+\sum\chi(C)\lesssim a_n\chi(\Ga)+a_n \sim a_n\chi(\Si)$.

\thm{Remarque}\fthm  On peut y ajouter un terme de ramification. La formule de Riemann-Hurewitz pour $f_n\vert_C$ donne $\chi(C)+\text{ram}(f_n\vert_C)=\text {deg}(f_n\vert_C)$. Soit $r_n=\sum \text{ram}(f_n\vert_C)$. C'est la ramification significative, celle loin du bord $\Bn$. Comme plus haut $\sum \text {deg}(f_n\vert_C)\leq a_n$. Donc $\chi(\Sn)+r_n\lesssim a_n\chi(\Si).$

\thm{Argument longueur-aire}\fthm

Rappelons pour finir comment une application enti\`ere non constante $f:\C\to \Sigma$ produit une suite de disques concentriques dans $\C$ satisfaisant l'hypoth\`ese. On note $h\vert dz \vert$ la pseudom\'etrique induite par $f$ sur $\C$, $a(r)$ l'aire du disque centr\'e en 0 de rayon $r$ et $l(r)$ la longueur de son bord. En coordonn\'ees polaires $a'=\int_0^{2\pi}h^2rd\theta$ et $l=\int_0^{2\pi}hrd\theta$. Par Cauchy-Schwarz $l^2\leq 2\pi r a'$. Donc $\int_1^{\infty}(\frac {l}{a})^2\frac {dr}{r} \leq \frac {2\pi}{a(1)}<+\infty$. On en d\'eduit bien une suite $r_n$ telle que $l(r_n)=o(a(r_n))$.

\Refs

\widestnumber\no{99}
\refno=0

\bref \by H. de Th\'elin\paper Une d\'emonstration du th\'eor\`eme de recouvrement des surfaces d'Ahlfors \jour Enseign.Math.\vol51\yr2005\pages203--209
\endref

\bref \by R. Nevanlinna \book Analytic functions \bookinfo Grund. Math. Wiss. \vol162 \publ Springer \yr 1970 \publaddr Berlin
\endref

\bref \by E. Reyssat \book Quelques aspects des surfaces de Riemann \bookinfo Progr. Math. \vol77 \publ Birkh\"auser \yr 1989 \publaddr Boston
\endref
 
\endRefs

\enddocument